\documentclass{article}

\usepackage[english]{babel}

\usepackage[letterpaper,top=2cm,bottom=2cm,left=3cm,right=3cm,marginparwidth=1.75cm]{geometry}

\usepackage{amsmath,amsthm,amssymb,mathrsfs,bm}
\usepackage{graphicx}
\usepackage[colorlinks=true, allcolors=blue]{hyperref}
\usepackage{doi}
\numberwithin{equation}{section}

\newtheorem{theorem}{\bf{Theorem}}[section]
\newtheorem{lemma}{\bf {Lemma}}[section]
\newtheorem{define}{\bf{Definition}}[section]

\newtheorem{proposition}{\bf{Proposition}}[section]
\newtheorem{remark}{\bf{Remark}}[section]

\newcommand{\be}{\begin{equation}}
\newcommand{\ee}{\end{equation}}
\newcommand\bes{\begin{eqnarray}}
\newcommand\ees{\end{eqnarray}}
\newcommand{\bess}{\begin{eqnarray*}}
\newcommand{\eess}{\end{eqnarray*}}
\newcommand{\mbE}{\hat{\mathbb{E}}}
\newcommand{\mbe}{\hat{\mathcal{E}}}
\newcommand{\V}{\mathbb{V}}
\newcommand{\mv}{\mathcal{V}}
\newcommand{\sles}{(\Omega,\mathcal{H},\mbE)}
\newcommand{\usigma}{\overline{\sigma}}
\newcommand{\lsigma}{\underline{\sigma}}
\newcommand*{\dif}{\mathop{}\!\mathrm{d}}
\newcommand{\Vstar}{\hat{\mathbb{V}}^*}
\newcommand{\mvstar}{\hat{\mathcal{V}}^*}

\title{Strassen's law of the  iterated logarithm under sub-linear expectations}
\author{Wang-yun Gu, Li-xin Zhang}
\date{}

\begin{document}
\maketitle

\begin{abstract}
\par We establish the Strassen's law of the iterated logarithm(LIL for short) for independent and identically distributed random variables with $\mbE[X_1]=\mbe[X_1]=0$ and $C_{\V}[X_1^2]<\infty$ under a sub-linear expectation space with a countably sub-additive capacity $\V$. We also show the LIL for upper capacity with $\sigma=\usigma$ under some certain conditions.\\
\par \bf{Keywords:}\rm\quad Sub-linear expectation, Capacity, Strassen's law of the iterated logarithm.
\end{abstract}

\section{Introduction}
Let $\{X_n;n\geq 1\}$ be a sequence of independent and identically distributed random variables on a probability space $(\Omega,\mathcal{F},P)$ with $E[X_1]=0,E[X_1^2]=1$. Denote $S_n=\sum_{k=1}^n X_k$ and 
\bess
\eta(t)=([t]+1-t)S_{[t]}+(t-[t])S_{[t]+1},
\eess
for any $t\geq 0$, which is the linear interpolation of $S_n$ at $n$. Strassen\cite{St64} established the classical  strong invariance principle for the law of iterated logarithm(LIL for short), which states that for $t\in[0,1]$ and $n\geq3$, with probability one,
\bess
\eta_n(t):=\frac{\eta(nt)}{(2n\log\log n)^{1/2}}
\eess
is relatively compact and the set of its limit points as $n$ tends to $\infty$ is K, where
\bess
K=\left\{f:f\in C[0,1]\enspace is\enspace absolutely\enspace continuous, f(0)=0,\int_0^1\vert f'(t)\vert^2dt\leq1\right\}.
\eess 
\par Numerous literatures extended Strassen's work since the classical independent and identically distributed condision may be too strong. Many papers have shown that it is sufficient for weak dependence, non-stationary or negatively dependence condition to establish the corresponding LIL.(see \cite{Mc77,Yu96,zhang01})
\par Recently, Wu and Chen\cite{Wu15} considered Strassen's LIL in a sub-linear expectation space $\sles$. Before their work, Hu and Chen\cite{Hu14} established the estimation of the increments of G-Brownian motion by means of  Cs\"{o}rg\H{o} and R\'{e}v\'{e}sz \cite{CR81} and obtained the LIL of G-Brownian motion as a corollary. Wu and Chen\cite{Wu15} proved the general LIL for G-Brownian motion. However, in \cite{Wu15}, they proved Strassen's LIL under the condition that $X_i$'s are bounded and continuous and did not give out either the lower bound of the clusters of $\eta_n$ or the LIL of upper capacity for $\eta_n(t)$. Moreover, the capacity in \cite{Wu15} is a special case of the form:
\bess
	\V^{\mathcal{P}}(A)=\sup_{P\in\mathcal{P}}P(A),\quad \forall A\in\mathcal{F},
\eess
where $\mathcal{P}$ is a weakly compact set of probability measures on $(\Omega,\mathcal{F})$.
\par The purpose of this paper is to establish Strassen's LIL under sub-linear expeactation space. In this paper, we will consider a countably sub-additive $\V$ and weaken the condition of $X_i$'s to 
\bess
C_{\V}[X_1^2]<\infty,
\eess
and we obtain both the upper and lower bound of the clusters of $\eta_n$ and the upper capacity form for $\sigma=\usigma$.
\par A difficulty in showing our theorem is that the truncation method has many limits. It will cause the mean uncertainty of the random variables. And if we take an uniform truncation $Y_i=(-c)\vee X_i\wedge c$, it is hard to control $\sum_{k=1}^n\mbE[X_k-Y_k]$. On the other hand, if we take $Y_i=(-c_i)\vee X_i\wedge c_i$, $Y_i$ will be in different distribution. So we have to make use of inequalities for sequences not necessarily identically distributed as far as possible. Fortunately, Zhang\cite{zhang21} has established these powerful inequalities for independent random variables under sub-linear expectations. Furthermore, Guo\cite{Guo22} and Zhang\cite{zhang21} enlightened us on the form of truncation. 
\par Zhang\cite{zhang21b} established a link between a sequence of independent random variables under the sub-linear expectation and a sequence of independent random variables on $\mathbb{R}^\infty$ under a probability measure $Q$. Based on this, he gave a purely probabilistic way to prove the conclusions under the upper cappacity. In the light of this idea, we derive the second part of our conclusion.
\par This paper is organized as follows. In Section 2, we introduce the basic notations of sub-linear expectation, capacity and Choquet integral. We recall some basic tools needed in our proof in Section 3. We also modify the exponential inequality for the maximum absolute value of independent random variables and obtain the increments of the linear interpolation of partial sums under sub-linear expectations. Two key propositions are obtained in Section 4 and Strassen's LIL under sub-linear expectations is established.

\section{Basic settings}
\par We use the framework and notations of Peng\cite{Peng08,Peng09,Peng19}.If one is familiar with these notations, he or she can skip this section. Let $(\Omega,\mathcal{F})$ be a given measurable space and let $\mathcal{H}$ be a linear space of real functions defined on  $(\Omega,\mathcal{F})$ such that if $X_1,\cdots, X_n\in\mathcal{H}$, then $\varphi(X_1,\cdots,X_n)\in\mathcal{H}$ for each $\varphi\in C_{l,Lip}(\mathbb{R}^n)$, where $C_{l,Lip}(\mathbb{R}^n)$ denotes the linear space of local Lipschitz functions $\varphi$ satisfying
\bess
\vert \varphi(\bm{x})-\varphi(\bm{y})\vert\leq C(1+\vert \bm{x}\vert^m+\vert\bm{y}\vert^m)\vert\bm{x}-\bm{y}\vert,\quad\forall\bm{x},\bm{y}\in\mathbb{R}^n,
\eess

~~~~~~~~~~~~~~~~~~~~~~~~~~~~~~~~~~for some $C>0, m\in\mathbb{N}$ depending on $\varphi$.\\
$\mathcal{H}$ is considered as a space of "random variables". We also denote $C_{b,Lip}(\mathbb{R}^n)$ the space of bounded Lipschitz functions. In this case, we denote $X\in\mathcal{H}$. 
\begin{define}
A sub-linear expectation $\mbE$ on $\mathcal{H}$ is a function $\mbE:\mathcal{H}\rightarrow\bar{\mathbb{R}}$ satisfying the following properties: for all $X,Y\in\mathcal{H}$, we have
\begin{itemize}
\item[(a)] Monotonicity: If $X\geq Y$, then $\mbe[X]\geq\mbe[Y]$;
\item[(b)] Constant preserving: $\mbE[c]=c$;
\item[(c)] Sub-additivity: $\mbE[X+Y]\leq\mbe[X]+\mbE[Y]$ whenever $\mbE[X]+\mbE[Y]$ is not of the form $+\infty-\infty$ or $-\infty+\infty$;
\item[(d)] Positive homogeneity: $\mbE[\lambda X]=\lambda\mbE[X]$.
\end{itemize}
Here, $\bar{\mathbb{R}}=[-\infty,\infty]$, $0\cdot\infty$ is defined to be 0. The triple $\sles$ is called a sub-linear expectation space. Give a sub-linear expectation $\mbE$, let us denote the conjugate expectation $\mbe$ of $\mbE$ by
\bess
\mbe[X]:=-\mbE[-X],\quad\forall X\in\mathcal{H}.
\eess
\end{define}
\par From the definition, it is easily shown that $\mbe[X]\leq\mbE[X],\enspace \mbE[X+c]=\mbE[X]+c$, and $\mbE[X-Y]\geq\mbE[X]-\mbE[Y]$ for all $X,Y\in\mathcal{H}$ with $\mbE[Y]$ being finite. We also call $\mbE[X]$ and $\mbe[X]$ the upper-expecation and lower-expectation of $X$, respectively.

\begin{define}
\begin{itemize}
\item[(i)] (Identical distribution) Let $\bm{X}_1$ and $\bm{X}_2$ be two n-dimensional random vectors, respectively, defined in sub-linear expectation spaces $(\Omega_1,\mathcal{H}_1,\mbE_1)$ and $(\Omega_2,\mathcal{H}_2,\mbE_2)$. They are called identically distributed, denoted by $\bm{X}_1\overset{d}{=}\bm{X}_2$, if
\bess
	\mbE_1[\varphi(\bm{X}_1)]=\mbE_2[\varphi(\bm{X}_2)],\quad\forall\varphi\in C_{l,Lip}(\mathbb{R}^n).
\eess
\par A sequence $\{X_n;n\geq 1\}$ of random variables is said to be identically distributed if $X_i\overset{d}{=}X_1$ for each $i\geq 1$.
\item[(ii)] (Independence) In a sub-linear expectation space $\sles$, a random vector $\bm{Y}=(Y_1,\cdots,Y_n),Y_i\in\mathcal{H}$ is said to be independent of another random vector $\bm{X}=(X_1,\cdots,X_m),X_i\in\mathcal{H}$ under $\mbE$ if for each test function $\varphi\in C_{l,Lip}(\mathbb{R}^m\times\mathbb{R}^n)$ we have $\mbE[\varphi(\bm{X},\bm{Y})]=\mbE[\mbE[\varphi(\bm{x},\bm{Y})]\vert_{\bm{x}=\bm{X}}]$, whenever $\bar{\varphi}(\bm{x}):=\mbE[\vert\varphi(\bm{x},\bm{Y})\vert]<\infty$ for all $\bm{x}$ and $\mbE[\vert\bar{\varphi}(\bm{X})\vert]<\infty.$
\item[(iii)] (Independent random variables) A sequence of random variables $\{X_n;n\geq 1\}$ is said to be independent if $X_{i+1}$ is independent of $(X_1,\cdots,X_i)$ for each $i\geq 1$.
\end{itemize}
\end{define}
It is easily seen that if $\{X_1,\cdots,X_n\}$ are independent, then $\mbE[\sum_{i=1}^nX_i]=\sum_{i=1}^nX_i$.
\par Next, we consider the capacities corresponding to the sub-linear expectations. Let $\mathcal{G}\subset\mathcal{F}$. A function $V:\mathcal{G}\rightarrow[0,1]$ is called a capacity if
$$
V(\emptyset)=0,\enspace V(\Omega)=1 \enspace and\enspace  V(A)\leq V(B)\enspace \forall A\subset B, A,B\in\mathcal{G}.
$$
It is called sub-additive if $V(A\cup B)\leq V(A)+V(B)$ for all $A,B\in\mathcal{G}$ with $A\cup B\in\mathcal{G}$.
\par Let $\sles$ be a sub-linear expectation space. We define $(\V,\mv)$ as a pair of capacities with the properties that
\bes
	\mbE[f]\leq\V(A)\leq\mbE[g]\quad if\enspace f\leq I_A\leq g,f,g\in\mathcal{H}\enspace and\enspace A\in\mathcal{F},\label{eq1.1}
\ees
$\V$ is sub-additive and $\mv(A):=1-\V(A^c),A\in\mathcal{F}$. It is obvious that
\bes
\mv(A\cup B)\leq\mv(A)+\V(B).
\ees
We call $\V$ and $\mv$ the upper and lower capacity, respectively. 
\par Also, we define the Choquet integrals/expectations $(C_{\V},C_{\mv})$ by
\bess
	C_V[X]=\int_0^{\infty}V(X\geq t)\dif t+\int_{-\infty}^0[V(X\geq t)-1]\dif t
\eess
with $V$ being replaced by $\V$ and $\mv$, respectively. It is obvious that $\mbE[X]\leq C_{\V}[X]$. If $\V$ on the sub-linear expectation space $\sles$ and $\tilde{\mathbb{V}}$ on the sub-linear expectation space $(\tilde{\Omega},\tilde{\mathcal{H}},\tilde{\mathbb{E}})$ are two capacities have the property (\ref{eq1.1}), then for any random variables $X\in\mathcal{H}$ and $\tilde{X}\in\tilde{\mathcal{H}}$ with $X\overset{d}{=}\tilde{X}$, we have
\bes
	\V(X\geq x+\epsilon)\leq\tilde{V}(\tilde{X}\geq x)\leq\V(X\geq x-\epsilon)\quad for\enspace all\enspace\epsilon>0\enspace and \enspace x.\label{eq1.2}
\ees
In fact, let $f\in C_{b,Lip}(\mathbb{R})$ such that $I\{y\geq x+\epsilon\}\leq f(y)\leq I\{y\geq x\}$. Then
\bess
	\V(X\geq x+\epsilon)\leq\mbE[f(X)]=\tilde{\mathbb{E}}[f(\tilde{X})]\leq\tilde{\V}(X\geq x),
\eess
and similar $\tilde{V}(\tilde{X}\geq x)\leq\V(X\geq x-\epsilon)$. From (\ref{eq1.2}), it follows that $\V(X\geq x)=\tilde{\V}(\tilde{X}\geq x)$ if $x$ is a continuous point of functions $\V(X\geq y)$ and $\tilde{\V}(\tilde{X}\geq y)$. Since a monotone function has at most countable number of discontinuous points, so
\bess
\V(X\geq x)=\tilde{\V}(\tilde{X}\geq x)
\eess 
for all but except countable many $x$, and then
\bess
C_{\V}[X]=C_{\tilde{\V}}[X].
\eess
\par In general, we choose $(\V,\mv)$ as
\bes
\hat{\V}(A):=\inf\{\mbE[\xi]:I_A\leq\xi,\xi\in\mathcal{H}\},\hat{\mv}(A)=1-\hat{\V}(A^c), \enspace\forall A\in\mathcal{F}.
\ees
\par Since $\hat{\V}$ may be not countably sub-additive so that the Borel-Cantelli lemma is not valid, we consider its countably sub-additive extension $\Vstar$ which defined by
\bes
\Vstar(A):=\inf\left\{\sum_{n=1}^{\infty}\hat{\V}(A_n):A\subset\bigcup_{n=1}^{\infty}A_n\right\},\mvstar(A)=1-\Vstar(A^c),\quad A\in\mathcal{F}.
\ees
As shown in Zhang\cite{zhang16}, $\Vstar$ is countably sub-additive, and $\Vstar(A)\leq\hat{\V}(A)$. Further, $\hat{\V}(A)$(resp. $\Vstar$) is the largest sub-additive(resp. countably sub-additive) capacity in sense that if $V$ is also a sub-additive(resp. countably sub-additive) capacity satisfying $V(A)\leq\mbE[g]$ whenenver $I_A\leq g\in\mathcal{H}$, then $V(A)\leq\hat{\V}(A)$(resp. $V(A)\leq\Vstar(A)$).
\par Besides $\hat{\V}^*$, another countably sub-additive capacity generated by $\mbE$ can be defined as follows:
\bes
	\mathbb{C}^*=\inf\left\{\lim_{n\rightarrow\infty}\mbE[\sum_{i=1}^n g_i]:I_A\leq\sum_{i=1}^\infty g_n,0\leq g_n\in\mathcal{H}\right\},A\in\mathcal{F}.
\ees
Then $\mathbb{C}^*\leq\hat{\V}^*$. It can be shown that the out capacity $c'$ defined in Example 6.5.1 of Peng\cite{Peng19} coincides with $\mathbb{C}^*$ if $\mathcal{H}$ is chosen as the family of (bounded) continuous functions on a metric space $\Omega$.
\par Futhermore, a random variable $X$ is called tight(under a capacity $\V$ satisfying (\ref{eq1.1})) if $\V(\vert X\vert\geq c)\rightarrow0$ as $c\rightarrow0$. It is obvious that if $\mbE[\vert X\vert]<\infty$ or $C_{\V}[\vert X\vert]<\infty$, then $X$ is tight.
\par Through this paper, for real numbers $x$ and $y$, we denote $\log(x)=\ln \max\{e,x\},x\vee y=\max\{x,y\},x\wedge y=\min\{x,y\},x^+=x\vee 0$ and $x^-=x\wedge 0$. For a random variable $X$, because $XI\{\vert X\vert\leq c\}$ may not be in $\mathcal{H}$, we will truncate it in the form $(-c)\vee X\wedge c$ denoted by $X^{(c)}$.

\section{Related inequalities and properties}
Exponential inequalities, Kolmogorov's converse exponential inequality and Borel-Cantelli lemma are basic tools established for the law of the iterated logarithm, the corresponding results under sub-linear expectation are obtained in Zhang\cite{zhang21}. We recall the lemmas needed in this paper and give a special case of exponential inequality.
\begin{lemma}
(\cite{zhang21}Lemma3.1) Let $\{X_1,\cdots,X_n\}$ be a sequence of independent random variables in the sub-linear expectation space $\sles$. Set $S_n=\sum_{i=1}^n X_i, B_n^2=\sum_{i=1}^n \mbE[X_i^2]$. Denote $B^2_{n,y}=\sum_{i=1}^n\mbE[(X_i\wedge y)^2], y>0.$ Then for all $x,y>0$,
\bes
	\V\left(\max_{k\leq n}(S_k-\mbE[S_k])\geq x\right)\leq\V\left(\max_{k\leq n}X_k>y\right)+\exp\left\{-\frac{x^2}{2(xy+B^2_{n,y})}\left(1+\frac{2}{3}\ln\left(1+\frac{xy}{B^2_{n,y}}\right)\right)\right\}.
\ees
\label{l4}
\end{lemma}

Next lemma is a special case for the absolute value of partial sums.
\begin{lemma}
Let $\{X_1,\cdots,X_n\}$ be a sequence of independent random variables in the sub-linear expectation space $\sles$. Set $S_n=\sum_{i=1}^n X_i, B_n^2=\sum_{i=1}^n \mbE[X_i^2]$. Then for all $x>(\max_{k\leq n}\mbE[S_k])^+\vee(\min_{k\leq n}\mbe[S_k])^-$and $y>0$,
\bes
	\V\left(\max_{k\leq n}\vert S_k\vert\geq x\right)&\leq& 2\V\left(\max_{k\leq n}\vert X_k\vert>y\right)+\exp\left\{-\frac{(x-\max_{k\leq n}\mbE[S_k])^2}{2(y(x-\max_{k\leq n}\mbE[S_k])+B^2_{n})}\right\}\nonumber\\
&+&\exp\left\{-\frac{(x+\min_{k\leq n}\mbe[S_k])^2}{2(y(x+\min_{k\leq n}\mbe[S_k])+B^2_{n})}\right\}.
\ees
\label{l1}
\end{lemma}
\proof Denote $B^2_{n,y}=\sum_{i=1}^n\mbE[(X_i\wedge y)^2], y>0$. Applying Lemma \ref{l4}, we have
\bess
	\V\left(\max_{k\leq n} S_k\geq x\right)&\leq& \V\left(\max_{ k\leq n}(S_k-\mbE[S_k]+\mbE[S_k])\geq x\right)\\
&\leq&\V\left(\max_{ k\leq n}(S_k-\mbE[S_k])\geq x-\max_{k\leq n}\mbE[S_k]\right)\\
&\leq&\V\left(\max_{k\leq n} X_k>y\right)+\exp\left\{-\frac{(x-\max_{k\leq n}\mbE[S_k])^2}{2(y(x-\max_{k\leq n}\mbE[S_k])+B^2_{n,y})}\right\}\\
&\leq& \V\left(\max_{k\leq n}\vert X_k\vert>y\right)+\exp\left\{-\frac{(x-\max_{k\leq n}\mbE[S_k])^2}{2(y(x-\max_{k\leq n}\mbE[S_k])+B^2_{n})}\right\}.
\eess
Similarly we can the symmetric inequality by replacing $X_i$ with $-X_i$, and thus the lemma is proved. \qedsymbol

\begin{remark}
	If we take $x=x_n$ such that $\lim_{n\rightarrow\infty}\frac{\sum_{k=1}^n\vert\mbE[X_k]\vert+\sum_{k=1}^n\vert\mbe[X_k]\vert}{x_n}=0$, when $n$ is large enough, it follows that
\bes
	\V\left(\max_{k\leq n} \vert S_k\vert\geq x\right)&\leq&2\V\left(\max_{k\leq n}\vert X_k\vert>y\right)+2\exp\left\{-\frac{x_n^2(1+o(1))}{2(x_ny+B^2_{n})}\right\}.
\ees

\end{remark}
The remark is obvious since 
\bess
	\vert\max_{k\leq n}\mbE[S_k]\vert\leq\sum_{k=1}^n\vert\mbE[X_k]\vert\enspace and\enspace \vert\min_{k\leq n}\mbe[S_k]\vert\leq\sum_{k=1}^n\vert\mbe[X_k]\vert.
\eess

The following lemma gives a lower bound of an exponential inequality under $\mv$ for independent and identically distributed random variables.
\begin{lemma}
	(\cite{zhang21}Lemma3.3) Suppose ${X_{ni};i=1,\cdots,k_n}$ is an array of independent and identically distributed random variables in the sub-linear expectation space $\sles$ with
\bess
\mbE[X^2_{n1}]\rightarrow \usigma^2<\infty, \mbe[X^2_{n1}]\rightarrow \lsigma^2>0,
\eess
and
\bess
\frac{\sum_{i=1}^{k_n}(\vert\mbE[X_{ni}]\vert+\vert\mbe[X_{ni}]\vert)}{x_n\sqrt{k_n}}=\frac{\sqrt{k_n}(\vert\mbE[X_{n1}]\vert+\vert\mbe[X_{n1}]\vert)}{x_n}\rightarrow 0.
\eess
Denote $S_n=\sum_{i=1}^{k_n}X_{ni}$. Then for any $z>0$,
\bes
	\liminf_{n\rightarrow\infty} x_n^{-2}\ln\mv\left(S_n\geq z\lsigma x_n\sqrt{k_n}\right)\geq -\frac{z^2}{2}.
\ees
\label{l3}
\end{lemma}

Borel-Cantelli lemma and its converse is still true under the upper and lower capacities $\V$ and $\mv$. 
\begin{lemma}
(\cite{zhang21} Lemma4.1)
\begin{itemize}
	\item[(1)] Let $\{A_n,n\geq 1\}$ be a sequence of events in $\mathcal{F}$. Suppose that $\V$ is a sub-additive capacity and $\sum_{n=1}^\infty \V(A_n)<\infty$. Then
	\bess
	\lim_{n\rightarrow\infty}\max_N\V\left(\bigcup_{i=n}^\infty A_n\right)=0.
	\eess
	If $\V$ is a countably sub-additive capacity, then
	\bess
	\V(A_n,i.o.)=0.
	\eess
	\item[(2)]Suppose $\{\xi_n;n\geq 1\}$ is a sequence of independent random variables in $\sles$ with a countably sub-additive $\V$. Then
	\bess
		\mv(\{\xi_n\geq 1\} i.o.)=1\enspace if\enspace  \sum_{n=1}^\infty \mv(\xi_n\geq 1+\epsilon)=\infty\enspace for\enspace some\enspace \epsilon>0.
	\eess
\end{itemize}
\end{lemma}
Next property is mentioned in Zhang\cite{zhang21} and is obvious by the definition of the lower capacity.
\begin{remark}
 If $\V$ is a countably sub-additive capacity, then the lower capacity $ \mv$ has the following property:
\bes
	\mv\left(\bigcap_{i=1}^\infty A_i\right)=1 \enspace for \enspace events \enspace \{A_n\}\enspace with\enspace A_n\supset A_{n+1}\enspace and\enspace \mv(A_n)=1, n=1,2,\cdots.\label{eq1}
\ees 
\end{remark}

The following lemma is useful to estimate the increments of the partial sums under lower capacity and the main idea of the proof is from Zhang\cite{zhang01}.
\begin{lemma}
	Let $\{X_n;n\geq 1\}$ be a sequence of independent random variables in the sub-linear expectation space $\sles$ with a countably sub-additive capacity $\V$. Suppose there exists a sequence of positive constants $\{\kappa_n;n\geq 1\}$ with $\kappa_n\downarrow 0$ and $\kappa_n\sqrt{\frac{n}{\log\log n}}\rightarrow\infty$ such that $\vert X_n\vert\leq \kappa_n\sqrt{\frac{n}{\log\log n}}$.\\
Furthermore, assume that $\mbE[X_1^2]\leq\usigma^2$ and $\lim_{n\rightarrow\infty}\frac{\sum_{k=1}^n\vert\mbE[X_k]\vert+\sum_{k=1}^n\vert\mbe[X_k]\vert}{(2n\log\log n)^{1/2}}=0$. Denote $T_n=\sum_{i=1}^n X_i$ and $B_n^2=\sum_{i=1}^n\mbE[X_i^2]$. Let $\{a_N;N\geq 1\}$ be a non-decreasing sequence of integers for which
\begin{itemize}
	\item[(1)]$1\leq a_N\leq N$,
	\item[(2)]$\lim_{N\rightarrow\infty}N/a_N<\infty$,
	\item[(3)]$N/a_N^\tau$ is eventually non-decreasing for some $\tau>0$.
\end{itemize} 
 Denote $\beta_N=\{2a_N(\log(N/a_N)+\log\log N)\}^{-\frac{1}{2}}$. Then we have
\bes
	\mv\left(\limsup_{N\rightarrow\infty}\max_{1\leq n\leq N}\max_{0\leq k\leq a_N}\beta_N\vert T_{n+k}-T_n\vert\leq 3\usigma\right)=1.\label{eq2}
\ees
\label{l2}
\end{lemma}

\proof
Let $x_N=(3/2+\epsilon/2)\usigma(2a_N(\log(N/a_N)+\log\log N))^{1/2}$, it follows that
\bess
\lim_{N\rightarrow\infty}\frac{\sum_{k=1}^N\vert\mbE[X_k]\vert+\sum_{k=1}^n\vert\mbe[X_k]\vert}{x_N}=\lim_{N\rightarrow\infty}\frac{\sum_{k=1}^N\vert\mbE[X_k]\vert+\sum_{k=1}^n\vert\mbe[X_k]\vert}{(2N\log\log N)^{1/2}}=0.
\eess

By applying Lemma \ref{l1}, we have, for any $\epsilon>0$, when $N$ is large enough,
\bess
	&&\V\left(\max_{1\leq n\leq N}\max_{0\leq k\leq a_N}\beta_N\vert T_{n+k}-T_n\vert\geq(3+\epsilon)\usigma\right)\\
&\leq&\V\left(\max_{0\leq l\leq[\frac{N}{a_N}]-1}\max_{la_N\leq n\leq(l+1)a_N}\max_{0\leq k\leq a_N}\beta_N\vert T_{n+k}-T_n\vert\geq(3+\epsilon)\usigma\right)\\
&\leq&\left[\frac{N}{a_N}\right]\max_l \V\left(\max_{la_N\leq n\leq(l+1)a_N}\max_{0\leq k\leq a_N}\beta_N\vert T_{n+k}-T_{la_N}-(T_n-T_{la_N})\vert\geq(3+\epsilon)\usigma \right)\\
&\leq&2\frac{N}{a_N}\max_l \V\left(\max_{0\leq k\leq 2a_N}\beta_N\vert T_{la_N+k}-T_la_N\vert\geq(3/2+\epsilon/2)\usigma \right)\\
&\leq&C\frac{N}{a_N}\exp \left\{-\frac{(1+o(1))(3/2+\epsilon/2)^2\usigma^2 a_N(\log(N/a_N)+\log\log N)}{(3/2+\epsilon/2)\kappa_{2a_N}\sqrt{\frac{2a_N}{\log\log2a_N}}\usigma/\beta_N+B_{2a_N}^2}\right\}\\
&\leq&C\frac{N}{a_N}\exp \left\{-\frac{(1+o(1))(3/2+\epsilon/2)^2\usigma^2(\log(N/a_N)+\log\log N)}{(3+\epsilon)\kappa_{2a_N}\sqrt{\frac{\log(N/a_N)+\log\log N}{\log\log2a_N}}\usigma+2\usigma^2}\right\}.
\eess
Note that $\kappa_{2a_N}\sqrt{\frac{\log(N/a_N)+\log\log N}{\log\log2a_N}}\rightarrow 0$ and $\frac{(3/2+\epsilon/2)^2}{2}\geq 1+\frac{\epsilon}{2}$, it follows that
\bess
	&&\V\left(\max_{1\leq n\leq N}\max_{0\leq k\leq a_N}\beta_N\vert T_{n+k}-T_n\vert\geq(3+\epsilon)\usigma\right)\\
&\leq& C\frac{N}{a_N}\exp \left\{-(1+o(1))(1+\epsilon/2)(\log(N/a_N)+\log\log N)\right\}\\
&\leq& C\frac{N}{a_N}\exp \left\{-(1+\epsilon/4)(\log(N/a_N)+\log\log N)\right\}\\
&=& C\left(\frac{N}{a_N}\right)^{-\epsilon/4}\left(\log N\right)^{-(1+\epsilon/4)}\leq C(\log N)^{1+\epsilon/4}.
\eess
Let $N_k=[\theta^k]$ for some $\theta>1$, by Borel-Cantelli lemma, we have
\bess
	\V\left(\limsup_{k\rightarrow\infty}\max_{1\leq n\leq N_k}\max_{0\leq i\leq a_{N_k}}\beta_{N_k}\vert T_{n+i}-T_n\vert\geq (3+\epsilon)\usigma\right)=0,
\eess
which means
\bess
	\mv\left(\limsup_{k\rightarrow\infty}\max_{1\leq n\leq N_k}\max_{0\leq i\leq a_{N_k}}\beta_{N_k}\vert T_{n+i}-T_n\vert\leq (3+\epsilon)\usigma\right)=1.
\eess
Since $\V$ is countably sub-additive, $\mv$ has property (\ref{eq1}) and thus
\bess
	\mv\left(\limsup_{k\rightarrow\infty}\max_{1\leq n\leq N_k}\max_{0\leq i\leq a_{N_k}}\beta_{N_k}\vert T_{n+i}-T_n\vert\leq 3\usigma\right)=1.
\eess
When $k$ is large enough, it follows that $a_{N_k}\geq \theta^{-1/\tau}a_{N_{k+1}}$, hence
\bess
	&&\limsup_{N\rightarrow\infty}\max_{1\leq n\leq N}\max_{0\leq k\leq a_N}\beta_N\vert T_{n+k}-T_n\vert\\
&\leq&\limsup_{k\rightarrow\infty}\max_{1\leq n\leq N_{k+1}}\max_{0\leq i\leq a_{N_{k+1}}}\frac{\vert T_{n+i}-T_n\vert}{(2a_{N_k}(\log(N_k/a_{N_{k+1}})+\log\log N_k))^{1/2}}\\
&\leq& \theta^{\frac{1}{2\tau}}\limsup_{k\rightarrow\infty}\max_{1\leq n\leq N_{k+1}}\max_{0\leq i\leq a_{N_{k+1}}}\beta_{N_{k+1}}\vert T_{n+i}-T_n\vert.
\eess
Thus
\bess
	&&\mv\left(\limsup_{N\rightarrow\infty}\max_{1\leq n\leq N}\max_{0\leq k\leq a_N}\beta_N\vert T_{n+k}-T_n\vert\leq 3\theta^{\frac{1}{2\tau}}\usigma\right)\\
	&\geq&	\mv\left(\limsup_{k\rightarrow\infty}\max_{1\leq n\leq N_k}\max_{0\leq i\leq a_{N_k}}\beta_{N_k}\vert T_{n+i}-T_n\vert\leq 3\usigma\right)=1.
\eess
By the arbitrariness of $\theta>1$ and  property (\ref{eq1}) again, we obtain euqation (\ref{eq2}).$\qedsymbol$
\par The following lemma is a special case of Lemma3.1 in Guo\cite{Guo22} for independent and identically distributed random variables with $C_{\V}[X^2_1]<\infty$. 
\begin{lemma}
	Let $\{X_n;n\geq 1\}$ be a sequence of independent and identically distributed random variables in the sub-linear expectation space $\sles$ with $C_{\V}[X_1^2]<\infty$. Then there exists a sequence of positive constants $\{\kappa_i;i\leq1\}$ with $\kappa_i\downarrow0$ and $\kappa_i\sqrt{i/\log\log i}\rightarrow\infty$ such that
\bes
	\sum_{i=16}^\infty\frac{\mbE\left[\left(\vert X_i\vert-\kappa_i\sqrt{i/\log\log i}\right)^+\right]}{\sqrt{2i\log\log i}}<\infty.
\ees
\label{l5}
\end{lemma}

\par Now we introduce the Proposition 2.1 established in Zhang\cite{zhang21b} as a basic tool to build a link between the capacity and a probability measure.
\begin{lemma}
	Let $\sles$ be a sub-linear expectation space with a capacity $\V$ satisfying (\ref{eq1.1}), and $\{X_n;n\geq1\}$ be a sequence of random variables in $\sles$. We can find a new sub-linear space $(\tilde{\Omega},\tilde{\mathcal{H}},\tilde{\mathbb{E}})$ defined on a metric space $\tilde{\Omega}=\mathbb{R}^\infty$, with a sequence $\{\tilde{X}_n;n\geq1\}$ of random variables and a set function $\tilde{V}:\tilde{F}\rightarrow[0,1]$ on it satisfying the following properties, where $\tilde{F}=\sigma(\tilde{\mathcal{H}})$:
\begin{itemize}
\item[(a)] $(X_1,X_2,\cdots,X_n)\overset{d}{=}(\tilde{X}_1,\tilde{X}_2,\cdots,\tilde{X}_n),n=1,2,\cdots,i,e.,$
\bess
	\tilde{\mathbb{E}}[\varphi(\tilde{X}_1,\tilde{X}_2,\cdots,\tilde{X}_n)]=\mbE[\varphi(X_1,X_2,\cdots,X_n)],\varphi\in C_{l,Lip}(\mathbb{R}^n),n\geq 1,
\eess
whenever the sub-linear expectation in the right hand is finite. In particular, if $\{X_n;n\geq1\}$ are independent under $\mbE$, then $\{\tilde{X}_n;n\geq1\}$ are independent under $\tilde{\mathbb{E}}$.
\item[(b)] Define
\bess
	\tilde{V}(A)=\tilde{\V}^{\tilde{\mathcal{P}}}(A)=\sup_{P\in\tilde{\mathcal{P}}}P(A),A\in\tilde{\mathcal{F}},
\eess
where $\tilde{\mathcal{P}}$ is the family of all probability measures $P$ on $(\tilde{\Omega},\tilde{\mathcal{F}})$ with the property
\bess
	P[\varphi]\leq\tilde{\mathbb{E}}\enspace for\enspace bounded\enspace \varphi\in\tilde{\mathcal{H}},
\eess
and $\tilde{V}\equiv0$ if $\tilde{\mathcal{P}}$ is empty. Then $\tilde{V}:\tilde{F}\rightarrow[0,1]$ is a countably sub-additive and nondecreasing function, and $\tilde{V}\leq\tilde{\mathbb{C}}^*\leq\tilde{\V}^*\leq\tilde{\V}$, where $\tilde{\mathbb{C}}^*,\tilde{\V}^*$ and $\tilde{\V}$ are defined on $(\tilde{\Omega},\tilde{\mathcal{H}},\tilde{\mathbb{E}})$ in the same way as $\mathbb{C}^*,\hat{\V}^*$ and $\hat{\V}$ on $\sles$, respectively.\\
Here and in the sequel, for a probability measure $P$ and a measurable function $X$, $PX$ is defined to be the expectation $\int X\dif P$.
\item[(c)]  If each $X_n$ is tight, then $\tilde{\mathcal{P}}$ is a weakly compact family of probability measures on the metric space $\tilde{\Omega}$,
\bess
	\tilde{\mathbb{E}}[\varphi]=\sup_{P\in\tilde{\mathcal{P}}}P[\varphi]\enspace for\enspace bounded\enspace \varphi\in\tilde{\mathcal{H}},
\eess
and $\tilde{V}$ is a countably sub-additive capacity with the property (\ref{eq1.1}),i.e.,
\bess
\tilde{\mathbb{E}}[f]\leq\tilde{V}(A)\leq\tilde{\mathbb{E}}[g]\quad if\enspace f\leq I_A\leq g,f,g\in\tilde{\mathcal{H}}\enspace and\enspace A\in\tilde{\mathcal{F}}.
\eess
\item[(d)] If $\{X_n;n\geq1\}$ are independent under $\mbE$ and each $X_n$ is tight, then for any sequence of vectors $\{\bm{\xi}_k=(X_{n_{k-1}+1},\cdots,X_{n_k});k\geq1\}$ and a sequence $\{E_k;k\geq 1\}$ of finite additive linear expectations on $\mathcal{H}_b=\{f\in\mathcal{H};f\enspace is \enspace bounded\}$ with $E_k\leq\mbE$, where $1=n_0<n_1<\cdots$, there exists a probability measure $Q$ on $\tilde{\Omega}$ such that $\{\tilde{\bm{\xi}}_k=(\tilde{X}_{n_{k-1}+1},\cdots,\tilde{X}_{n_k});k\geq1\}$ is a sequence of independent random vectors under $Q$,
\bes
	Q[\varphi(\tilde{\bm{\xi}}_k)]=E_k[\varphi(\bm{\xi}_k)]\enspace for\enspace all \enspace \varphi\in C_{b,Lip}(\mathbb{R}^{n_k-n_{k-1}}),\\
Q[\varphi(\tilde{X_1},\cdots,\tilde{X_n})]\leq\mbE[\varphi(X_1,\cdots,X_n)] ]\enspace for\enspace all \enspace \varphi\in C_{b,Lip}(\mathbb{R}^{n}),
\ees
and
\bes
\tilde{v}\left((\tilde{X}_1,\tilde{X}_2,\cdots)\in B\right)\leq Q\left((\tilde{X}_1,\tilde{X}_2,\cdots)\in B\right)\leq\tilde{V}\left((\tilde{X}_1,\tilde{X}_2,\cdots)\in B\right)
\ees
for all $B\in\mathcal{B}(\mathbb{R}^\infty)$, where $\tilde{v}(A)=1-\tilde{V}(A^c)$.
\end{itemize}
\label{l6}
\end{lemma}

The following lemma is useful to estimate the difference between the original and truncated random variables. 
\begin{lemma}
Suppose $X\in\mathcal{H},1\leq p\leq2,C_{\V}(\vert X\vert^p)<\infty.$ Then
\bes
\mbE[(\vert X\vert-c)^+]=o(c^{1-p}).
\ees
\label{l7}
\end{lemma}
The proof is similar to Lemma 4.1 of Zhang\cite{zhang21b}.

Next lemma can be seen as an extension of Kronecker lemma.
\begin{lemma}
Let $\{x_n;n\geq1\}$ be a sequence of real number with $x_n\rightarrow x$ as $n\rightarrow\infty$, for a fixed positive integer $d$, $\alpha_1,\cdots,\alpha_d$ are real numbers such that $\alpha_1^2+\cdots+\alpha_d^2=1$, then it follows that 
\bes
\frac{d}{N}\sum_{i=1}^d\alpha_i^2\sum_{j=[(i-1)N/d]+1}^{[iN/d]}x_j\rightarrow x \label{eq3.1}
\ees
as $N\rightarrow\infty$.\label{l8}
\end{lemma}
\proof For any $\epsilon>0$, there exists $N_0\in\mathbb{N}$ such that for all $n>N_0$, we have $\vert x_n-x\vert<\epsilon$, and it is obvious that 
\bess
x=\lim_{N\rightarrow \infty}\frac{d}{N}\sum_{i=1}^d\alpha_i^2\sum_{j=[(i-1)N/d]+1}^{[iN/d]}x.
\eess
Hence we have
\bess
&&\lim_{N\rightarrow \infty}\left\vert \frac{d}{N}\sum_{i=1}^d\alpha_i^2\sum_{j=[(i-1)N/d]+1}^{[iN/d]}x_j-x\right\vert\\
&\leq&\lim_{N\rightarrow \infty}\left\vert \frac{d}{N}\sum_{i=1}^d\alpha_i^2\sum_{j=[(i-1)N/d]+1}^{[iN/d]}(x_j-x)\right\vert+\lim_{N\rightarrow \infty}\left\vert \frac{d}{N}\sum_{i=1}^d\alpha_i^2\sum_{j=[(i-1)N/d]+1}^{[iN/d]}x-x\right\vert\\
&\leq&\lim_{N\rightarrow \infty} \frac{d}{N}\sum_{i=1}^d\alpha_i^2\sum_{j=[(i-1)N/d]+1}^{[iN/d]}\vert x_j-x\vert\\
&=&\lim_{N\rightarrow \infty} \frac{d}{N}\sum_{i=1}^d\alpha_i^2\sum_{\overset{j=[(i-1)N/d]+1}{j>N_0}}^{[iN/d]}\vert x_j-x\vert\\
&\leq&\lim_{N\rightarrow \infty} \frac{d}{N}\sum_{i=1}^d\alpha_i^2\sum_{j=[(i-1)N/d]+1}^{[iN/d]}\epsilon=\epsilon.
\eess
By the arbitrariness of $\epsilon$, we have (\ref{eq3.1}).\qedsymbol

The following lemma shows that under certain conditions, the inverse Borel-Cantelli lemma for some special $\V$ still holds without the assumption of continuity.
\begin{lemma}
(Zhang\cite{zhang21b} Lemma 2.5)Let $\sles$ be a sub-linear expectation space with a capacity $V$ having the property (\ref{eq1.1}), and $v(A)=1-V(A^c)$. Suppose one of the following conditions is satisfied.
\begin{itemize}
	\item[(a)]The sub-linear expectation $\mbE$ satisfies
	\bess
		\mbE[X]=\max_{P\in\mathcal{P}}P[X],\enspace X\in\mathcal{H}_b,
	\eess
		where $\mathcal{H}_b=\{f\in\mathcal{H}:f\enspace is\enspace bounded\}$, $\mathcal{P}$ is a countable-dimensionally weakly compact family of probability measures on $(\Omega,\sigma(\mathcal{H}))$ in sense that, for any $Y_1, Y_2,\cdots\in\mathcal{H}_b$ and any sequence $\{P_n\}\subset\mathcal{P}$, there is a subsequence $\{n_k\}$ and a probability measure $P\in\mathcal{P}$ for which
	\bess
		\lim_{k\rightarrow\infty}P_{n_k}[\varphi(Y_1,\cdots,Y_d)]=P[\varphi(Y_1,\cdots,Y_d)],\enspace \varphi\in C_{b,Lip}(\mathbb{R}^d),d\geq1.
	\eess
	\item[(b)] $\mbE$ on $\mathcal{H}_b$ is regular in sense that $\mbE[X_n]\downarrow0$ for any elements $\mathcal{H}_b\ni X_n\downarrow0$. Let $\mathcal{P}$ be the family of all probability measures on $(\Omega,\sigma(\mathcal{H}))$ for which
	\bess
		P[f]\leq\mbE[f],f\in\mathcal{H}_b.
	\eess
	\item[(c)] $\Omega$ is a complete separable metric space, each element $X(\omega)$ in $\mathcal{H}$ is a continuous function on $\Omega$. The capacity $V$ with the property (\ref{eq1.1}) is tight in sense that, for any $\epsilon>0$, there is a compact set $K\subset\Omega$ such that $V(K^c)<\epsilon$. Let $\mathcal{P}$ be defined as in (b).
	\item[(d)] $\Omega$ is a complete separable metric space, each element $X(\omega)$ in $\mathcal{H}$ is a continuous function on $\Omega$. The sub-linear expectation $\mbE$ is defined by
	\bess
		\mbE[X]=\max_{P\in\mathcal{P}}P[X],
	\eess
where $\mathcal{P}$ is a weakly compact family of probability measures on the metric space $\Omega$.
\end{itemize}
Denote $\V^{\mathcal{P}}(A)=\max_{P\in\mathcal{P}}P(A),A\in\sigma(\mathcal{H})$. Let $\{\bm{X}_n;n\geq 1\}$ be a sequence of independent random vectors in $\sles$, where $\bm{X}_n$ is $d_n$-dimensional. If $F_n$ is a $d_n-dimensional$ close set with $\sum_{n=1}^\infty v(\bm{X}_n\notin F_n)<\infty,$ then for $\V=\V^{\mathcal{P}},\mathbb{C}^*,\hat{\V}^*$ or $\hat{\V}$,
\bess
	\mv(\bm{X}_n\notin F_{n}\enspace i.o.)=0;
\eess
If $F_{n,j}$s are $d_n$-dimensional close sets with $\sum_{n=1}^\infty V(\bm{X}_n\in F_{n,j})=\infty, j=1,2,\cdots$, then for $\V=\V^{\mathcal{P}},\mathbb{C}^*,\hat{\V}^*$ or $\hat{\V}$,
\bess
	\V\left(\bigcap_{j=1}^\infty \{\bm{X}_n\in F_{n,j} \enspace i.o.\}\right)=1.
\eess

\label{l9}
\end{lemma}

\section{Main results}
\par In this section, we will consider the Strassen's invariance principle for LIL under sub-linear expectation space $\sles$ with a countably sub-additive capacity $\V$. Let $C[0,1]$ denote the space of all continuous maps from $[0,1]$ to $\mathbb{R}$. For any $\sigma>0$, denote
\bess
	K_\sigma=\left\{f:f\in C[0,1]\enspace is\enspace absolutely\enspace continuous, f(0)=0,\int_0^1\vert f'(t)\vert^2dt\leq\sigma^2\right\}.
\eess
\begin{theorem}
	Let $\{X_n;n\geq 1\}$ be independent and identically distributed random variables in the sub-linear expectation space $\sles$ with $\mbE[X_1]=\mbe[X_1]=0,C_{\V}[X_1^2]<\infty$ and $\V$ is a countably sub-additive capacity on $\sles$. Denote $\mbE[X_1^2]=\usigma^2, \mbe[X_1^2]=\lsigma^2$ and $S_n=\sum_{k=1}^n X_k$, and for $t\in[0,1]$
	\begin{align*}
		\begin{split}
 			\eta_n(t)= \left \{
 			\begin{array}{ll}
 				(2n\log\log n)^{-1/2}S_k   & t=\frac{k}{n},\\
    			linear,     & in\enspace between.
			\end{array}
		\right.
		\end{split}
	\end{align*}
Then $\{\eta_n(t);n\geq3\}$is relatively compact and
\bes
	\mv\left( K_{\lsigma}\subset C\left\{\eta_n(t)\right\}\subset K_{\usigma} \right)=1,\label{eq4.6}
\ees
Further, if the space $\sles$ satisfies one of the conditions (a)-(d) in Lemma \ref{l9}, then for $\V=\V^{\mathcal{P}},\mathbb{C}^*$ or $\hat{\V}^*$
\begin{align}
		\begin{split}
 			\V\left(C\left\{\eta_n(t)\right\}=K_\sigma\right)= \left \{
 			\begin{array}{ll}
 				1,     & \sigma=\usigma,\\
    			0,     & \sigma\notin [\lsigma,\usigma].\label{eq4.7}
			\end{array}
		\right.
		\end{split}
\end{align}
where $C\{x_n\}$ denotes the cluster set of a sequence of $\{x_n\}$ in $\mathbb{R}$.\label{th1}
\end{theorem}

The following propositions are of great significance to prove Theorem \ref{th1}.
\begin{proposition}
	Let $\{X_n;n\geq 1\}$ be a sequence of independent and identically distributed random variables in the sub-linear expectation space $\sles$ with $\mbE[X_1]=\mbe[X_1]=0,C_{\V}[X_1^2]<\infty$, and $\V$ is a countably sub-additive capacity on $\sles$. Denote $  \mbE[X_1^2]=\usigma^2,\mbe[X_1^2]=\lsigma^2.$ Let $d\geq1$ be a positive integer and $\alpha_1,\cdots,\alpha_d$ be real numbers with $\alpha_1^2+\cdots+\alpha_d^2=1$. Denote $S_n=\sum_{k=1}^n X_k$ and
\bess
T_N=\sum_{i=1}^d \alpha_i(S_{[iN/d]}-S_{[(i-1)N/d]}).
\eess
Then
\bes
	\mv\left(\lsigma\leq\limsup_{N\rightarrow\infty}\frac{\sqrt{d}T_N}{(2N\log\log N)^{1/2}}\leq\usigma\right)=1
	\label{eq3}
\ees
and
\bes
	\mv\left(-\usigma\leq\liminf_{N\rightarrow\infty}\frac{\sqrt{d}T_N}{(2N\log\log N)^{1/2}}\leq-\lsigma\right)=1.
\label{eq4}
\ees
\label{p1}
\end{proposition}
\proof We only give the proof of (\ref{eq3}), since (\ref{eq4}) follows from (\ref{eq3}) immediately by symmetry.  Without loss of generality, we can assume that $\lsigma>0$. Otherwise, the left part of the conclusion is trivial by  noting (\ref{eq1}) and for all $m$,
\bes
\mv\left(\bigcup_{N=m}^\infty \frac{\sqrt{d}T_N}{(2N\log\log N)^{1/2}}\geq-\epsilon\right)&\geq&\mv\left(\frac{\sqrt{d}T_N}{(2N\log\log N)^{1/2}}\geq-\epsilon\right)\nonumber\\
&=&1-\V\left(\frac{-\sqrt{d}T_N}{(2N\log\log N)^{1/2}}>\epsilon\right)\nonumber\\
&\geq&1-\frac{d\mbE[T_N^2]}{2\epsilon^2 N\log\log N}\rightarrow 1.\label{eq4.1}
\ees

\par In order to prove (\ref{eq3}), it is sufficient to show that for any $0<\epsilon<1/4$,
\bes
	\mv\left(\limsup_{N\rightarrow\infty}\frac{\sqrt{d}\vert T_N\vert}{(2N\log\log N)^{1/2}}\leq(1+4\epsilon)\usigma\right)=1
\label{eq5}
\ees
and
\bes
	\mv\left(\limsup_{N\rightarrow\infty}\frac{\sqrt{d} T_N}{(2N\log\log N)^{1/2}}\geq(1-4\epsilon)\lsigma\right)=1.
	\label{eq6}
\ees
\par First, we show (\ref{eq5}). Let $c_n=\kappa_n\sqrt{\frac{n}{\log\log n}}$ and $Y_n=(-c_n)\vee X_n\wedge c_n$, where $\kappa_n$ is the positive constants constructed in Lemma \ref{l5}. By Lemma \ref{l5} and Kronecker lemma, we obtain
\bes
	\lim_{n\rightarrow\infty}\frac{\sum_{k=1}^n \mbE[\vert X_k-Y_k\vert]}{\sqrt{2n\log\log n}}=0,\label{eq12}
\ees
By (\ref{eq12}) and similar to (\ref{eq4.1}), we have
\bes
	\mv\left(\lim_{n\rightarrow\infty}\frac{\sum_{k=1}^n\vert X_k-Y_k\vert}{\sqrt{2n\log\log n}}=0\right)=1.\label{eq8}
\ees
Since $\mbE[X_k]=\mbe[X_k]=0$, we have
\bes
	\vert\max_{k\leq n}\sum_{i=1}^k\mbE[Y_i]\vert+\vert\min_{k\leq n}\sum_{i=1}^k\mbe[Y_i]\vert&\leq&\sum_{k=1}^n\vert\mbE[Y_k]\vert+\sum_{k=1}^n\vert\mbe[Y_k]\vert\nonumber\\
&=&\sum_{k=1}^n\vert\mbE[Y_k-X_k]\vert+\sum_{k=1}^n\vert\mbe[Y_k-X_k]\vert \nonumber\\
&\leq&2\sum_{k=1}^n\mbE[\vert Y_k-X_k\vert]=o(\sqrt{2n\log\log n}).\label{eq13}
\ees

Define
\bess
\bar{S}_n&=&\sum_{k=1}^n Y_k,\quad \hat{S}_n=\sum_{k=1}^n (X_k-Y_k),\\
\bar{T}_N&=&\sum_{i=1}^d \alpha_i(\bar{S}_{[iN/d]}-\bar{S}_{[(i-1)N/d]}),\\
\hat{T}_N&=&\sum_{i=1}^d \alpha_i(\hat{S}_{[iN/d]}-\hat{S}_{[(i-1)N/d]}).
\eess

By Cauchy-Schwarz inequality, we have
\bes
	\limsup_{N\rightarrow\infty}\frac{\sqrt{d}\vert\hat{T}_N\vert}{(2N\log\log N)^{1/2}}&\leq&\limsup_{N\rightarrow\infty}\frac{\sqrt{d}(\sum_{i=1}^d(\hat{S}_{[iN/d]}-\hat{S}_{[(i-1)N/d]})^2)^{1/2}}{(2N\log\log N)^{1/2}}\nonumber\\
&\leq&\limsup_{N\rightarrow\infty}\frac{d\max_{1\leq i\leq d}\vert\hat{S}_{[iN/d]}-\hat{S}_{[(i-1)N/d]}\vert}{(2N\log\log N)^{1/2}}\nonumber\\
&\leq&2d\limsup_{N\rightarrow\infty}\frac{\max_{1\leq i\leq d}\vert\hat{S}_{[iN/d]}\vert}{(2N\log\log N)^{1/2}}\nonumber\\
&\leq&2d\limsup_{N\rightarrow\infty}\frac{\sum_{i=1}^d\vert\hat{S}_{[iN/d]}\vert}{(2N\log\log N)^{1/2}}\nonumber\\
&\leq&2d^2\limsup_{N\rightarrow\infty}\frac{\vert\hat{S}_N\vert}{(2N\log\log N)^{1/2}}.\label{eq9}
\ees
It follows from (\ref{eq8}) and (\ref{eq9}) that
\bes
	\mv\left(\lim_{N\rightarrow\infty}\frac{\sqrt{d}\vert\hat{T}_N\vert}{\sqrt{2N\log\log N}}=0\right)=1.\label{eq4.2}
\ees
We only need to show that
\bes
	\mv\left(\limsup_{N\rightarrow\infty}\frac{\sqrt{d}\vert \bar{T}_N\vert}{(2N\log\log N)^{1/2}}\leq(1+3\epsilon)\usigma\right)=1
\label{eq10}
\ees

We still set $N_k=[\theta^k]$ for some $\theta>1$. From Lemma \ref{l2}, we have 
\bess
	&&\limsup_{k\rightarrow\infty}\max_{dN_k\leq N\leq dN_{k+1}}\frac{\sqrt{d}\vert \bar{T}_N-\bar{T}_{dN_k}\vert}{(2N\log\log N)^{1/2}}\\
&\leq&\limsup_{k\rightarrow\infty}\max_{dN_k\leq N\leq dN_{k+1}}\frac{\vert \bar{T}_N-\bar{T}_{dN_k}\vert}{(2N_k\log\log N_k)^{1/2}}\\
&\leq&\limsup_{k\rightarrow\infty}\max_{0\leq n\leq N_{k+1}}\max_{0\leq m\leq N_{k+1}-N_k}\frac{\vert \bar{S}_{n+m}-\bar{S}_{n}\vert}{(2N_k\log\log N_k)^{1/2}}\\
&\leq&3\usigma(\theta-1)^{\frac12}\rightarrow 0
\eess
as $\theta\rightarrow1$.\\

\par Thus we only need to show
\bes
\sum_{k=1}^\infty\V(\vert \bar{T}_{dN_k}\vert\geq(1+\epsilon)\usigma(2N_k\log\log N_k)^{1/2})<\infty. 
\label{eq7}
\ees

\par Since $\bar{T}_{dN_k}=\sum_{l=1}^d\sum_{i=(l-1)N_k+1}^{lN_k}\alpha_l Y_i$ and $\mbE[Y_i^2]\leq\usigma^2$, we denote $B_k^2=\sum_{l=1}^d\sum_{i=(l-1)N_k+1}^{lN_k}\mbE[\alpha_l^2 Y_i^2]\leq N_k\usigma^2$. Note $\kappa_{i}\rightarrow0$ as $i\rightarrow\infty$, so for any given $\epsilon>0$, there exists some $\alpha\leq\frac{2\epsilon+\epsilon^2}{2(1+\epsilon)}\usigma$, when $k$ is large enough, we have $\vert \alpha_lY_i\vert\leq \vert Y_i\vert\leq\kappa_i(\frac{i}{\log\log i})^{\frac12}\leq\alpha(\frac{N_k}{\log\log N_k})^{\frac12}$ for all $l=1,\cdots,d,i=1,\cdots,dN_k$.

\par Then by Lemma \ref{l1} with $x=(1+\epsilon)\usigma(2N_k\log\log N_k)^{1/2}$ and $y=\alpha(\frac{N_k}{\log\log N_k})^{\frac12}$, it follows that
\bess
	&&\V(\vert \bar{T}_{dN_k}\vert\geq(1+\epsilon)\usigma(2N_k\log\log N_k)^{1/2})\\
	&\leq&C\exp\left\{-\frac{(1+o(1))(1+\epsilon)^2\usigma\log\log N_k}{\sqrt{2}(1+\epsilon)\alpha+\usigma}\right\}\\
	&\leq&C\exp\left\{-(1+o(1))(1+\frac{(\sqrt{2}-1)(2\epsilon+\epsilon^2)}{\sqrt{2}+2\epsilon+\epsilon^2})\log\log N_k\right\}\\
	&\leq&C\exp\left\{-(1+\frac{2\epsilon+\epsilon^2}{3(\sqrt{2}+2\epsilon+\epsilon^2)})\log\log N_k\right\}.
\eess
So Equation (\ref{eq7}) is proved.

\par On the other hand, to prove (\ref{eq6}), we need another truncation. Let $N_k=[\theta^k]$ for some $\theta>d/\epsilon$ and $I(k)=\{N_k+1,\cdots,N_{k+1}\}$. For $j\in I(k)$, let $Z_j=(-c_{N_{k+1}})\vee X_j\wedge c_{N_{k+1}}$. Define
\bess
	\tilde{S}_n=\sum_{k=1}^n Z_k,\quad\check{S}_n=\sum_{k=1}^n(X_k-Z_k),\\
\tilde{T}_N=\sum_{i=1}^d \alpha_i(\tilde{S}_{[iN/d]}-\tilde{S}_{[(i-1)N/d]}),\\
\check{T}_N=\sum_{i=1}^d \alpha_i(\check{S}_{[iN/d]}-\check{S}_{[(i-1)N/d]}).
\eess
Note $\vert Z_j-X_j\vert\leq\vert Y_j-X_j\vert$, it follows from (\ref{eq12}),(\ref{eq8}) and the proof of (\ref{eq13}),(\ref{eq9}) that
\bess
		\lim_{n\rightarrow\infty}\frac{\sum_{j=1}^n\vert\mbE[Z_j]\vert+\sum_{j=1}^n\vert\mbe[Z_j]\vert}{\sqrt{2n\log\log n}}=0,\\
		\mv\left(\lim_{n\rightarrow\infty}\frac{\sqrt{d}\vert\check{T}_N\vert}{\sqrt{2n\log\log n}}=0\right)=1.
\eess

So it is sufficient to show that
\bess
	\mv\left(\limsup_{N\rightarrow\infty}\frac{\sqrt{d} \tilde{T}_N}{(2N\log\log N)^{1/2}}\geq(1-3\epsilon)\lsigma\right)=1.
\eess
Note
\bess
\limsup_{k\rightarrow\infty}\frac{\tilde{T}_{dN_k}}{(2N_k\log\log N_k)^{1/2}}=\limsup_{k\rightarrow\infty}\frac{\sqrt{d}\tilde{T}_{dN_k}}{(2dN_k\log\log dN_k)^{1/2}}\leq\limsup_{N\rightarrow\infty}\frac{\sqrt{d}\tilde{T}_{N}}{(2N\log\log N)^{1/2}},
\eess
we only need to show
\bes
	\mv\left(\limsup_{k\rightarrow\infty}\frac{ \tilde{T}_{dN_k}}{(2N_k\log\log N_k)^{1/2}}\geq(1-2\epsilon-\sqrt{\epsilon})\lsigma\right)=1,
	\label{eq11}
\ees
which is equivalent to 
\bess
\mv\left(\limsup_{k\rightarrow\infty}\frac{ \tilde{T}_{dN_k}}{(2N_k\log\log N_k)^{1/2}}\geq(1-3\epsilon)\lsigma\right)=1,
\eess

Let $z=1-2\epsilon, x_n=\sqrt{\log\log n}$ and
\bess
	W_k=\tilde{T}_{dN_k}-\alpha_1\tilde{S}_{dN_{k-1}}=\alpha_1(\tilde{S}_{N_k}-\tilde{S}_{dN_{k-1}})+\sum_{l=2}^d\alpha_l(\tilde{S}_{lN_k}-\tilde{S}_{(l-1)N_{k}}),
\eess
then we have
\bess
\left\{\frac{W_k}{\lsigma\sqrt{2N_k\log\log N_k}}\geq z\right\}&=&\left\{\frac{\alpha_1(\tilde{S}_{N_k}-\tilde{S}_{dN_{k-1}})+\sum_{l=2}^d \alpha_l(\tilde{S}_{lN_k}-\tilde{S}_{(l-1)N_k})}{\lsigma\sqrt{2N_k\log\log N_k}}\geq z\right\}\\
&\supset&\left\{ \frac{\alpha_1(\tilde{S}_{N_k}-\tilde{S}_{dN_{k-1}})}{\lsigma\sqrt{2N_k\log\log N_k}}\geq \alpha_1^2z\right\}\bigcap_{l=2}^d\left\{\frac{\alpha_l(\tilde{S}_{lN_k}-\tilde{S}_{(l-1)N_k})}{\lsigma\sqrt{2N_k\log\log N_k}}\geq \alpha_l^2 z\right\}.
\eess

\par For given $z>0$ and $\epsilon>0$, let $f\in C_{b,Lip}(\mathbb{R})$ such that $I\{x\geq z+\epsilon\}\leq f(x)\leq I\{x\geq z\}$. It follows that
\bess
	I\left\{\frac{W_k}{\lsigma\sqrt{2N_k\log\log N_k}}\geq z\right\}\geq f\left(\frac{\tilde{S}_{N_k}-\tilde{S}_{dN_{k-1}}}{\alpha_1\lsigma\sqrt{2N_k\log\log N_k}}\right)\prod_{l=2}^d f\left(\frac{\tilde{S}_{lN_k}-\tilde{S}_{(l-1)N_k}}{\alpha_l\lsigma\sqrt{2N_k\log\log N_k}}\right).
\eess

\par Note that $\{\tilde{S}_{N_k}-\tilde{S}_{dN_{k-1}}\},\{\tilde{S}_{lN_k}-\tilde{S}_{(l-1)N_k},l=2,\cdots,d\}$ are independent under $\mbE$ and $\mbe$, we can obtain
\bess
\mv\left(\frac{W_k}{\lsigma\sqrt{2N_k\log\log N_k}}\geq z\right)&\geq&\mbe\left[f\left(\frac{\tilde{S}_{N_k}-\tilde{S}_{dN_{k-1}}}{\alpha_1\lsigma\sqrt{2N_k\log\log N_k}}\right)\right]\prod_{l=2}^d \mbe\left[f\left(\frac{\tilde{S}_{lN_k}-\tilde{S}_{(l-1)N_k}}{\alpha_l\lsigma\sqrt{2N_k\log\log N_k}}\right)\right]\\
&\geq&\mv\left(\frac{\tilde{S}_{N_k}-\tilde{S}_{dN_{k-1}}}{\alpha_1\lsigma\sqrt{2N_k\log\log N_k}}\geq z+\epsilon\right)\prod_{l=2}^d \mv\left(\frac{\tilde{S}_{lN_k}-\tilde{S}_{(l-1)N_k}}{\alpha_l\lsigma\sqrt{2N_k\log\log N_k}}\geq z+\epsilon\right).
\eess
Without loss of generality, we can assume that $\alpha_l>0$, otherwise, we can remove $\alpha_l$ for $\alpha_l=0$ and replace $X_i$ with $-X_i$ if $\alpha_l<0$. Applying Lemma \ref{l3}, where the conditions are obviously satisfied, it follows that
\bess
	&&\liminf_{k\rightarrow\infty}x_{N_k}^{-2}\ln\mv\left(W_k\geq z\lsigma\sqrt{2N_k\log\log N_k}\right)\\
&\geq& \liminf_{k\rightarrow\infty}x_{N_k}^{-2} \{\ln\mv\left(\tilde{S}_{N_k}-\tilde{S}_{dN_{k-1}}\geq(1-\epsilon)\alpha_1\lsigma\sqrt{2N_k\log\log N_k}\right)\\
&+&\sum_{l=2}^d\ln\mv\left(\tilde{S}_{lN_k}-\tilde{S}_{(l-1)N_k}\geq(1-\epsilon)\alpha_l\lsigma\sqrt{2N_k\log\log N_k}\right)\}\\
&\geq&\liminf_{k\rightarrow\infty}x_{N_k}^{-2}\ln\mv\left(\tilde{S}_{N_k}-\tilde{S}_{dN_{k-1}}\geq(\sqrt{\frac{2\theta}{\theta-d}}(1-\epsilon)\alpha_1)\lsigma\sqrt{N_k-dN_{k-1}}x_{N_k}\right) \\
&+&\sum_{l=2}^d \liminf_{k\rightarrow\infty}x_{N_k}^{-2}\ln\mv\left(\tilde{S}_{lN_k}-\tilde{S}_{(l-1)N_k}\geq (\sqrt{2}(1-\epsilon)\alpha_l)\lsigma\sqrt{N_k}x_{N_k}\right)\\
&\geq&-(1-\epsilon)^2\alpha_1^2(\frac{\theta}{\theta-d})+\sum_{l=2}^d -(1-\epsilon)^2\alpha_l^2\\
&\geq&-(1-\epsilon)^2(\frac{\theta}{\theta-d})>-(1-\epsilon).
\eess
Therefore,
\bes
\sum_{k=1}^\infty \mv\left(W_k\geq(1-2\epsilon)\lsigma(2N_k\log\log N_k)^{1/2}\right)=\infty,
\label{eq14}
\ees
By (\ref{eq14}) and Borel-Cantelli lemma, we have
\bes
	\mv\left(\limsup_{k\rightarrow\infty}\frac{W_k}{\sqrt{2N_k\log\log N_k}}\geq(1-2\epsilon)\lsigma\right)=1.
\label{eq15}
\ees
\par For $N_k\leq n< N_{k+1}$, it follows that
\bes
	\limsup_{k\rightarrow\infty}\frac{\tilde{T}_{dN_k}}{\sqrt{2N_k\log\log N_k}}\geq \limsup_{k\rightarrow\infty}\frac{W_k}{\sqrt{2N_k\log\log N_k}}-\limsup_{k\rightarrow\infty}\frac{\vert\tilde{S}_{dN_{k-1}}\vert}{\sqrt{2N_k\log\log N_k}},
\label{eq16}
\ees
and
\bes
	\limsup_{k\rightarrow\infty}\frac{\vert\tilde{S}_{dN_{k-1}}\vert}{\sqrt{2N_k\log\log N_k}}=\limsup_{k\rightarrow\infty}\frac{\vert S_{dN_{k-1}}\vert}{\sqrt{2N_k\log\log N_k}}\leq\sqrt{\frac{d}{\theta}}\usigma<\sqrt{\epsilon}\usigma\rightarrow 0,
\label{eq17}
\ees
where the inequality follows from the law of iterated logarithm for independent and identically distributed random variables.
\par From (\ref{eq15}), (\ref{eq16}) and (\ref{eq17}), we can obtain that (\ref{eq11}) holds.$\qedsymbol$

\begin{proposition}
Under the conditions and notations in Proposition \ref{p1}, and we further assume that $\sles$ satisfies one of conditions (a)-(d) in Lemma \ref{l9}. Then for $\V=\V^{\mathcal{P}},\mathbb{C}^*$ or $\hat{\V}^*$, we have 
\bes
\V\left(\limsup_{N\rightarrow\infty}\frac{\sqrt{d}T_N}{(2N\log\log N)^{1/2}}=\usigma \right)=1.
\ees
\label{p2}
\end{proposition}
\proof Note (\ref{eq4.2}), it suffices to show 
\bes
\V^{\mathcal{P}}\left(\limsup_{N\rightarrow\infty}\frac{\sqrt{d}\bar{T}_N}{(2N\log\log N)^{1/2}}=\usigma \right)=1.
\label{eq4.5}
\ees
Let $\sigma\in[\lsigma,\usigma]$, there exists $\alpha\in[0,1]$ such that $\sigma^2=\alpha\lsigma^2+(1-\alpha)\usigma^2$. We set $\sigma_i^2=\alpha\mbE[Y_i^2]+(1-\alpha)\mbe[Y_i^2]$. Then by Theorem 1.2.1 of Peng\cite{Peng19},  for each i, there exists $\theta_{i,1},\theta_{i,2}\in\Theta$ such that 
\bess
	E_{\theta_{i,1}}[Y_i^2]=\mbE[Y_i^2], E_{\theta_{i,2}}[Y_i^2]=\mbe[Y_i^2].
\eess
Define the linear operator $E_i=\alpha E_{\theta_i,1}+(1-\alpha)E_{\theta_i,2}$. Then 
\bess
	E_i[Y_i^2]=\sigma_i^2,E_i\leq\mbE.
\eess
Since each $Y_n$ is tight, By Lemma \ref{l6}, we can construct a copy $\{\tilde{Y}_n;n\geq1\}$ on $(\tilde{\Omega},\tilde{\mathcal{H}},\tilde{\mathbb{E}})$ of $\{Y_n;n\geq1\}$, a countably sub-additive capacity $\tilde{V}$ with property (\ref{eq1.1}) and a probability measure $Q$ on $\tilde{\Omega}$ such that $\{\tilde{Y}_n;n\geq1\}$ is a sequence of independent random variables under $Q$, and we have
\bess
	Q[\varphi(\tilde{Y}_i)]=E_i[\varphi(Y_i)]\enspace for\enspace all \enspace \varphi\in C_{b,Lip}(\mathbb{R}),\\
Q[\varphi(\tilde{Y_1},\cdots,\tilde{Y_n})]\leq\mbE[\varphi(Y_1,\cdots,Y_n)] ]\enspace for\enspace all \enspace \varphi\in C_{b,Lip}(\mathbb{R}^{n}),
\eess
and
\bes
\tilde{v}\left(B\right)\leq Q\left(B\right)\leq\tilde{V}\left(B\right) \enspace for\enspace all\enspace B\in\sigma(\tilde{Y}_1,\tilde{Y}_2,\cdots).\label{eq4.4}
\ees
Note $\vert E_i[(Y_i^{(c)})^2]-E_i[Y_i^2]\vert\leq\mbE[(Y_i^2-c^2)^+]\rightarrow0$ as $c\rightarrow0$, we have
\bes
Q[\tilde{Y}_i]&=&\lim_{c\rightarrow\infty}Q[\tilde{Y}_i^{(c)}]=\lim_{c\rightarrow\infty}E_i[Y_i^{(c)}]\leq\mbE[Y_i],\\
Q[\tilde{Y}_i^2]&=&\lim_{c\rightarrow\infty}Q[\tilde{Y}_i^2\wedge c^2]=\lim_{c\rightarrow\infty}E_i[Y_i^2\wedge c^2]=E_i[Y_i^2]=\sigma_i^2.
\ees
By the linearity of $E_i$, we can have
\bess
	\mbe[Y_i]\leq Q[\tilde{Y_i}]\leq\mbE[Y_i].
\eess
Note by Lemma \ref{l7} that 
\bess
	&&\vert \mbE[Y_i]\vert=\vert\mbE[Y_i-X_i]\vert\leq\mbE[\vert X_i-Y_i\vert]=\mbE[(\vert X_1\vert-\kappa_i\sqrt{\frac{i}{\log\log i}})^+]=o((\kappa_i\sqrt{\frac{i}{\log\log i}})^{-1})\rightarrow0,\\
&&\vert \mbE[Y_i^2]-\mbE[X_i^2]\vert\leq\mbE[(X_1^2-(\kappa_i\sqrt{\frac{i}{\log\log i}})^2)^+]=o(1)\rightarrow0.
\eess
Similarly, we have $\vert\mbe[Y_i]\vert\rightarrow0$ and $\vert \mbe[Y_i^2]-\mbe[X_i^2]\vert\rightarrow0$, and thus
\bes
	Q[\tilde{Y_i}]\rightarrow0,\quad Q[\tilde{Y}_i^2]=\sigma_i^2\rightarrow\sigma^2.\label{eq4.3}
\ees
We denote $\breve{Y}_i=\tilde{Y}_i-Q[\tilde{Y}_i],\enspace s_n^2=\sum_{i=1}^nQ[\breve{Y}_i^2]$ and
\bess
\tilde{S}_n=\sum_{k=1}^n \tilde{Y}_k,\quad\breve{S}_n=\sum_{k=1}^n\breve{Y}_k,\\
\tilde{T}_N=\sum_{i=1}^d \alpha_i(\tilde{S}_{[iN/d]}-\tilde{S}_{[(i-1)N/d]}),\\
\breve{T}_N=\sum_{i=1}^d \alpha_i(\breve{S}_{[iN/d]}-\breve{S}_{[(i-1)N/d]}).
\eess
We first assume that $\lsigma>0$ and thus $\sigma>0$. From (\ref{eq4.3}), it is obvious that
\bes
s_n^2\rightarrow\infty, \quad \frac{s_n^2}{n}\rightarrow\sigma<\infty,\quad i.e.\enspace s_n^2=O(n).
\ees
Then there exists a constant $C$ such that $C\sqrt{\frac{s_n^2}{\log\log s_n^2}}\geq\sqrt{\frac{n}{\log\log n}}$. It follows from (\ref{eq1.2}), (\ref{eq4.4}) and (\ref{eq4.3}) that
\bess
Q\left(\vert \breve{Y}_n\vert\leq2C\kappa_n\sqrt{\frac{s_n^2}{\log\log s_n^2}}\right)&\geq& Q\left(\vert \tilde{Y}_n\vert\leq2\kappa_n\sqrt{\frac{n}{\log\log n}}-\vert Q[\tilde{Y}_n]\vert\right)\\
&\geq&\tilde{v}\left(\vert \tilde{Y}_n\vert\leq2\kappa_n\sqrt{\frac{n}{\log\log n}}-\vert Q[\tilde{Y}_n]\vert\right)\\
&\geq&\mv\left(\vert Y_n\vert\leq\kappa_n\sqrt{\frac{n}{\log\log n}}\right)=1.
\eess 
Hence it follows from (\ref{eq4.3}), Kolmogorov's bounded LIL and Lemma \ref{l8} that
\bess
\limsup_{N\rightarrow\infty}\frac{\sqrt{d}\tilde{T}_N}{(2N\log\log N)^{1/2}}&=&\limsup_{N\rightarrow\infty}\frac{\sqrt{d}\breve{T}_N}{(2N\log\log N)^{1/2}}\\
&=&\lim_{N\rightarrow\infty}\left(\sum_{i=1}^d \frac{d}{N}\alpha_i^2(s_{[iN/d]}^2-s_{[(i-1)N/d]}^2)\right)^{1/2}\\
&=&\lim_{n\rightarrow\infty}\left(Q[\tilde{Y}_n^2]-Q[\tilde{Y}_n]^2\right)^{1/2}=\sigma\quad Q-a.s..
\eess

If $\sigma=0$, we can set $\breve{Y}_i=\tilde{Y}_i-Q[\tilde{Y}_i]+\epsilon\xi_i$ instead, with $\epsilon>0$ and $\{\xi_i;i\geq1\}$ be a sequence of independent standard normal random variables under $Q$ which is also independent to $\{\tilde{Y}_i;i\geq1\}$. Similarly, we have
\bess
	\limsup_{N\rightarrow\infty}\frac{\sqrt{d}\breve{T}_N}{(2N\log\log N)^{1/2}}=\epsilon\quad Q-a.s.,
\eess 
and thus
\bess
	\limsup_{N\rightarrow\infty}\frac{\sqrt{d}\tilde{T}_N}{(2N\log\log N)^{1/2}}=0\quad Q-a.s.
\eess 
by the arbitrariness of $\epsilon$.

Then for any $\epsilon>0,$ there exists a subsequence $\{N_k\}$ satisfying $\frac{N_{k+1}}{N_k}>d/\epsilon$, such that
\bes
Q\left(\lim_{k\rightarrow\infty}\frac{\tilde{T}_{dN_k}}{(2N_k\log\log N_k)^{1/2}}=\sigma\right)=1.
\ees 
Denote
\bess
\tilde{W}_k=\tilde{T}_{dN_k}-\alpha_1\tilde{S}_{dN_{k-1}},\bar{W}_k=\bar{T}_{dN_k}-\alpha_1\bar{S}_{dN_{k-1}}.
\eess
Similar to (\ref{eq17}), we have
\bess
\limsup_{k\rightarrow\infty}\frac{\tilde{W}_k}{(2N_k\log\log N_k)^{1/2}}&\geq&\lim_{k\rightarrow\infty}\frac{\tilde{T}_{dN_k}}{(2N_k\log\log N_k)^{1/2}}-\limsup_{k\rightarrow\infty}\frac{\vert\tilde{S}_{dN_{k-1}}\vert}{(2N_k\log\log N_k)^{1/2}}\nonumber\\
&\geq&(1-\sqrt{\epsilon})\sigma,\enspace Q-a.s.
\eess
Then by the definition of upper limit and Borel-Cantelli lemma, we have for any $\epsilon>0$,
\bess
	Q\left(\frac{\tilde{W}_k}{(2N_k\log\log N_k)^{1/2}}\geq(1-\sqrt{\epsilon})\sigma\enspace i.o.\right)=1,
\eess
and thus
\bess
	\sum_{k=1}^\infty Q\left(\frac{\tilde{W}_k}{(2N_k\log\log N_k)^{1/2}}\geq(1-\sqrt{\epsilon})\sigma\right)=\infty.
\eess
Hence 
\bess
	\sum_{k=1}^\infty \V^{\mathcal{P}}\left(\frac{\bar{W}_k}{(2N_k\log\log N_k)^{1/2}}\geq(1-2\sqrt{\epsilon})\sigma\right)&\geq& \sum_{k=1}^\infty \tilde{V}\left(\frac{\tilde{W}_k}{(2N_k\log\log N_k)^{1/2}}\geq(1-\sqrt{\epsilon})\sigma\right)\\
&\geq& \sum_{k=1}^\infty Q\left(\frac{\tilde{W}_k}{(2N_k\log\log N_k)^{1/2}}\geq(1-\sqrt{\epsilon})\sigma\right)=\infty.
\eess
Noting the independence of $\bar{W}_k$, for any sequence $\epsilon_j\downarrow 0$, by Lemma \ref{l9}, we have
\bess
\V^{\mathcal{P}}\left(\bigcap_{j=1}^\infty \{\frac{\bar{W}_k}{(2N_k\log\log N_k)^{1/2}}\geq(1-2\sqrt{\epsilon_j})\sigma \enspace i.o.\}\right)=1.
\eess
Since
\bess
	\left\{\bigcap_{j=1}^\infty \{\frac{\bar{W}_k}{(2N_k\log\log N_k)^{1/2}}\geq(1-2\sqrt{\epsilon_j})\sigma \enspace i.o.\}\right\}\subset\left\{\limsup_{k\rightarrow\infty}\frac{\bar{W}_k}{(2N_k\log\log N_k)^{1/2}}\geq\sigma \right\},
\eess
it follows that
\bes
\V^{\mathcal{P}}\left(\limsup_{k\rightarrow\infty}\frac{\bar{W}_k}{(2N_k\log\log N_k)^{1/2}}\geq\sigma\right)=1.\label{eq4.8}
\ees
By (\ref{eq17}) and property (\ref{eq1}) again, we obtain
\bes
	\hat{\mv}^*\left(\limsup_{k\rightarrow\infty}\frac{\vert\bar{S}_{dN_{k-1}}\vert}{(2N_k\log\log N_k)^{1/2}}=0\right)=1.\label{eq4.9}
\ees
It follows from (\ref{eq4.8}) and (\ref{eq4.9}) that 
\bess
\V^{\mathcal{P}}\left(\frac{\sqrt{d}\bar{T}_N}{(2N\log\log N)^{1/2}}\geq\sigma\right)\geq\V^{\mathcal{P}}\left(\frac{\bar{T}_{dN_k}}{(2N_k\log\log N_k)^{1/2}}\geq\sigma\right)=1.
\eess

Note that $\hat{\mv}^*\left(\limsup_{N\rightarrow\infty}\frac{\sqrt{d}\bar{T}_N}{(2N\log\log N)^{1/2}}\leq\usigma\right)=1$, we have for all $\sigma\in[\lsigma,\usigma]$,
\bess
\V^{\mathcal{P}}\left(\sigma\leq\limsup_{N\rightarrow\infty}\frac{\sqrt{d}\bar{T}_N}{(2N\log\log N)^{1/2}}\leq\usigma\right)=1.
\eess
Hence (\ref{eq4.5}) is proved.\qedsymbol

\par With the preparation above, we now start to prove Theorem \ref{th1}.\\
\bf{Proof of Theorem 4.1}\rm.  For any $f\in C[0,1]$, define
\begin{align*}
	\begin{split}
 		f^{(d)}(x)= \left \{
 		\begin{array}{ll}
 			f(i/d)   & x=\frac{i}{d},\\
    		linear,     & in \enspace between.
		\end{array}
	\right.
	\end{split}
\end{align*}
Let
\bess
	C_d&=&\{f^{(d)}:f\in C[0,1]\}\subset C[0,1],\\
	K_\sigma^d&=&\{f^{(d)}:f\in K_\sigma\}.
\eess
Then similar to Cs\"{o}rg\H{o} and R\'{e}v\'{e}sz \cite{CR81}, Proposition \ref{p1} and \ref{p2} are equivalent to the sequence $\{\eta_n^{(d)};n\geq 1\}$ which is relatively compact in $C_d$ quasi surely and 
\bess
	\mv\left(K_{\lsigma}^d\subset C\{\eta_n^{(d)}\}\subset K_{\usigma}^d \right)=1,\\
	\V\left(C\{\eta_n^{(d)}\}=K_{\usigma}^d\right)=1.
\eess

On the other hand, by Lemma \ref{l2} and (\ref{eq8}), 
\bess
	&&\limsup_{n\rightarrow\infty}\sup_{0\leq x\leq1}\vert \eta_n(x)-\eta_n^{(d)}(x)\vert\\
	&\leq&\limsup_{n\rightarrow\infty}\frac{\max_{0\leq i\leq n}\max_{0\leq k\leq[n/d]}\vert S_{i+k}-S_i\vert}{(2n\log\log n)^{1/2}}\\
	&\leq&\limsup_{n\rightarrow\infty}\frac{\max_{0\leq i\leq n}\max_{0\leq k\leq[n/d]}\vert \bar{S}_{i+k}-\bar{S}_i\vert}{(2n\log\log n)^{1/2}}+2\limsup_{n\rightarrow\infty}\frac{\max_{0\leq i\leq 2n}\vert \hat{S}_i\vert}{(2n\log\log n)^{1/2}}\\
	&\leq& 3d^{-1/2}\usigma+2\limsup_{n\rightarrow\infty}\frac{\sum_{i=1}^{2n}\vert X_i-Y_i \vert}{(2n\log\log n)^{1/2}}\rightarrow 0,\quad as\enspace d\rightarrow\infty,
\eess
where $Y_i,\bar{S}_i$ and $\hat{S}_i$ are defined as in the proof of Proposition \ref{p1}. Hence ,we obtain  (\ref{eq4.6}) and $\V=\V^{\mathcal{P}},\mathbb{C}^*$ or $\hat{\V}^*$
\bess
	\V(C\{\eta_n(t)\}=K_{\usigma})=1.
\eess
For $\sigma<\lsigma$, it is obvious that for all countably sub-additive capacity $\V$,
\bess
	\V(C\{\eta_n(t)\}=K_\sigma)\leq\V(C\{\eta_n(t)\}\subset K_{\lsigma})\leq 1-\mv(C\{\eta_n(t)\}\subset K_{\lsigma})=0,
\eess
and for $\sigma>\usigma$ is similar. Hence (\ref{eq4.7}) holds for $\hat{\V}^*$ and thus $\V^{\mathcal{P}}$ and $\mathbb{C}^*$. The proof is completed.\qedsymbol

\bibliographystyle{plain}
\bibliography{SSIP}

\begin{thebibliography}{10}
\bibitem{CR81}
Cs{\"{o}}rg\H{o} M. and R\'{e}v\'{e}sz P.
\newblock {\em Strong Approximation in Probability and Statistics.}
\newblock Academic Press, New York, 1981.

\bibitem{Guo22}
Guo~X. F., Li~S., and Li~X. P.
\newblock On the hartman-wintner law of the iterated logarithm under sublinear
  expectation.
\newblock {\em Communications in Statistics-Theory and Methods}, 2022.
\newblock \doi{10.1080/03610926.2022.2026394}.

\bibitem{Hu14}
Hu~F., Chen~Z. J., and Zhang~D. F.
\newblock How big are the increments of g-brownian motion?
\newblock {\em Science in China-Mathematics}, 57(8):1687--1700, 2014.



\bibitem{Mc77}
McLeish~D. L.
\newblock On the invariance principle for nonstationary mixingales.
\newblock {\em The Annals of Probability}, 5(4):616--621, 1977.

\bibitem{Peng08}
Peng~S. G.
\newblock A new central limit throrem under sublinear expectations.
\newblock arXiv:0803.2656v1, 2008.

\bibitem{Peng09}
Peng~S. G.
\newblock Survey on normal distributions, central limit throrem, brownian
  motion and the related stochastic calculus under sublinear expectations.
\newblock {\em Sci. China Ser. A}, 52(7):1391--1411, 2009.

\bibitem{Peng19}
Peng~S. G.
\newblock {\em Nonlinear Expectations and Stochastic Calculus under Uncertainty
  with Robust CLT and G-Brownian Motion}.
\newblock Springer, Berlin, Heidelberg, 2019.







\bibitem{St64}
Strassen V.
\newblock An invariance principle for the law of the iterated logarithm.
\newblock {\em Z. Wahrscheinlichkeitstheorie}, 3:211--226, 1964.

\bibitem{Wu15}
Wu~P. Y. and Chen~Z. J.
\newblock Invariance principles for the law of the iterated logarithm under
  g-framework.
\newblock {\em Science in China-Mathematics}, 58(6):1251--1264, 2015.

\bibitem{Yu96}
Yu~H.
\newblock A strong invariance principle for associated sequences.
\newblock {\em The Annals of Probability}, 24(4):2079--2097, 1996.

\bibitem{zhang01}
Zhang~L. X.
\newblock Strassen's law of the iterated logarithm for negatively associated
  random vectors.
\newblock {\em Stochastic Process and their Applications}, 95:311--328, 2001.

\bibitem{zhang16}
Zhang~L. X.
\newblock Rosenthal's inequalities for independent and negatively dependent
  random variables under sub-linear expectations with applications.
\newblock {\em Science in China-Mathematics}, 59(4):751--768, 2016.

\bibitem{zhang21}
Zhang~L. X.
\newblock On the laws of the iterated logarithm under sub-linear expectations.
\newblock {\em Probability, Uncertainty and Quantitative Risk}, 6(4):409--460,
  2021.

\bibitem{zhang21b}
Zhang~L. X.
\newblock The sufficient and necessary conditions of the strong law of large numbers under the sub-linear expectations.
\newblock arXiv:2103.01390v3, 2021.



\end{thebibliography}

\end{document}